\documentclass[12pt]{amsart}
\usepackage[T2A,T1]{fontenc}
\usepackage[utf8]{inputenc}
\usepackage[russian,english]{babel}
\usepackage{amsmath,amssymb}
\usepackage[new]{old-arrows}
\usepackage{graphicx}
\usepackage[colorlinks,allcolors=blue]{hyperref}
\usepackage{tikz-cd} 
\usepackage{wrapfig}
\usepackage{caption}
\usepackage{lipsum}

\input epsf

\makeatletter
\newcommand{\insertheaderrule}{\rlap{\rule[-.3\normalbaselineskip]{\textwidth}{.4pt}}}
\let\old@evenhead\@evenhead \let\old@oddhead\@oddhead
\def\@evenhead{\insertheaderrule\old@evenhead}
\def\@oddhead{\insertheaderrule\old@oddhead}
\makeatother

\def\Isom{\mathop{\text{\rm Isom}}}
\def\T{\mathop{\text{\rm T}}}
\def\SO{\mathop{\text{\rm SO}}}
\def\H{\mathop{\text{\rm H}}}
\def\Hom{\mathop{\text{\rm Hom}}}
\def\Span{\mathop{\text{\rm Span}}}
\def\so3{\mathfrak{so}_3}

\begin{document}

\selectlanguage{english}

\markleft{\hskip43pt\hfill version 1.12\hfill02.10.2021}
\markright{09.04.2022\hfill version 1.12\hfill\hskip42pt}

\quad

\vskip-50pt

\centerline{\Large\bf Moduli spaces of polygons and}
\centerline{\Large\bf deformations of polyhedra with boundary}

\vskip20pt

\centerline{\large Sasha Anan$'$in$^\dagger$, Dmitrii Korshunov\footnote{\href{mailto:dmitrii.korshunov@impa.br}{dmitrii.korshunov@impa.br}} }

\vskip40pt

{\hfill\it To Fedor Bogomolov}

{\hfill\it with all our admiration}

{\footnotesize 
\begin{center}
\vspace{0.5cm}
    \textbf{Abstract}
\end{center}

We prove a conjecture of Ian Agol: all isometric realizations of a polyhedral surface with boundary sweep out an isotropic subset in the Kapovich-Millson moduli space of polygons isomorphic to the boundary. For a generic polyhedral disk we show that boundaries of its isometric realizations make up a Lagrangian subset. As an application of this result, we obtain a new solution to the problem of Richard Kenyon about spanning domes of piecewise linear curves comprised of unit intervals in $\mathbb R^3$ .
}

\vskip30pt

\centerline{\bf1.~Introduction}

\medskip

The space of isometric maps (not necessarily embeddings or immersions) of a triangulated metric polyhedron to $\mathbb{R}^3$ up to the action of the group of isometries of $\mathbb R^3$ is a scheme over the real numbers, given by a set of quadratic algebraic equations. There exists an extensive body of literature devoted to the study of the structure of this scheme. Even the question of its dimension (``rigidity and flexibility'') for polyhedra homeomorphic to closed surfaces is highly nontrivial. For an introduction to this subject see I.Kh. Sabitov's survey \cite{Sabitov} and references therein.

On the other hand, the space of spatial polygons up to isometries of $\mathbb R^3$ carries a lot of structure: it is a complex analytic space and its smooth locus is a symplectic manifold (Deligne-Mostow \cite{DM}, Kapovich-Millson \cite{KaM}, Klyachko \cite{Klyachko}). Singularities, if present, are isolated and quadratic \cite{KaM, Klyachko}. 

The boundary of a polyhedral surface is a polygon. Thus, there is a natural map $\delta$ from ``the moduli space of polyhedra'' to the ``moduli space of polygons'', which sends a polyhedron to its boundary. This is an algebraic map between two schemes over $\mathbb R$. However, the notion of Zariski tangent space (``the space of infinitesimal deformations'') allows us to speak about the derivative of $\delta$, skew-symmetric non-degenerate forms and pullbacks of forms even at singular points.

 The aim of this paper is to study the symplectic geometric properties of $\delta$. In particular, we show that the pull-back  of the Kapovich-Millson symplectic form with respect to $\delta$ is null (``$\delta$ is isotropic''). Moreover, for a generic polyhedral disk $P$, the image of $\delta$ is Lagrangian subset in the moduli space of polygons (that is, contains an open dense Lagrangian submanifold). We apply these results, following a suggestion of Ian Agol, to resolve a problem of Richard Kenyon about spanning domes of integral curves \cite{Kenyon}.

The basic idea behind the proof is to extend the notion of a polyhedron to the one that we call a {\it graph-surface}. Graph-surfaces have a well-defined notion of boundary and the map $\delta$. Their main advantage is that they are a common generalization of polyhedral surfaces with connected boundary and graphs. The class of graph-surfaces, unlike the class of polyhedral surfaces, is closed with respect to the operation of collapse described in the next paragraph. 

The second main ingredient is the notion of a {\it collapse} of a face adjacent to the boundary. It mutates one graph-surface $S$ to another one $S'$ with one face less. Let us denote maps
between corresponding spaces of polyhedra and polygons by $\delta_S$ and $\delta_{S'}$ respectively. Using the fact that a triangle is rigid, we will show that a collapse gives an embedding of the space of infinitesimal deformations of $S$ to that of $S'$ that induces a symplectic embedding between ranges of $d\delta_{S}$ and $d\delta_{S'}$ in the corresponding tangent spaces to the respective spaces of polygons. Thus we will prove that if the pull-back of $\delta_{S'}$ is null then the pull-back of $\delta_S$ is also null. 
Some sequence of collapses eventually terminates with a graph, that is, a graph-surface without faces. At last, the pull-back of $\delta$ for graphs can be obtained by a direct computation, and this completes the induction. 

The statement about the co-isotropic property of $\delta(S)$ in generic case  is proved by dimension counting. Wherever we use the notions of surfaces, triangulations, tubular neighborhoods, etc., PL-category is assumed.

The authors would like to thank Misha Verbitsky for his interest in this work and essential help during all stages of the preparation of the paper. The second author (D.K.) wants to express his gratitude to Konstantin Loginov for the argument in the proof of Lemma 3.2.2. The first author (S.A.) passed away during the preparation of the manuscript.

\bigskip

\centerline{\bf2.~The space $\mathbb E^P$ of polygons and the space $\mathbb E^S$ of polyhedra}

\medskip

In this section we introduce the main objects of our interest --- the space $\mathbb E^P$ of polygons, the space $\mathbb E^S$
of polyhedra, and closely related spaces. All polygons and polyhedra live in the $3$-dimensional euclidean space
$\mathbb E:=\mathbb E^3$.

\medskip

{\bf2.1~Graph-surface.} Let $\hat S$ be a closed surface and let $S\subset \hat S$ be a finite two-dimensional simplicial complex
with nondegenerate triangles and edges such that $D:=\hat S\smallsetminus S$ is homeomorphic to an open disk. The triangulation consists of the
set $V$ of vertices, of the set $E$ of edges taken with all orientations and equipped with the function $-:E\to E$, $e\mapsto-e$,
that flips the orientation, and of the set $T$ of triangles. Such an $S$ is a {\it graph-surface.}

\medskip

\noindent
\includegraphics[width=0.24\textwidth]{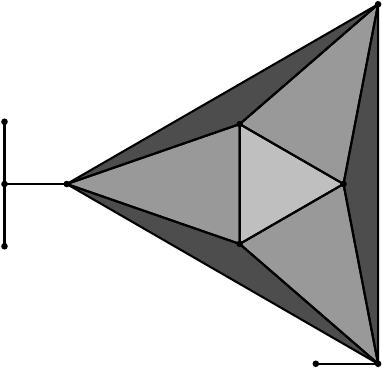}

\leftskip125pt

\vskip-114pt

A graph-surface in the simplest case is a usual triangulated subsurface in $\hat S$ bounding exactly one disk in $\hat S$. What can
happen in the general case is sketched on the picture on the left (where $\hat S$ is a $2$-sphere).

\medskip

{\bf2.2.~Lemma.} {\sl In the settings of\/ {\rm2.1}, the boundary\/ $\partial D$ of\/ $D$ is a union of edges. Moreover,
$\partial D$ can be decomposed into the union $\partial D=g_1\cup\dots\cup g_k$, with a prescribed cyclic order of the edges
$g_1,\dots,g_k\in E${\rm;} in this order any next edge begins at the end of the previous one. The list $g_1,\dots,g_k$ of edges
admits repetitions, with the\break}

\leftskip0pt

\vskip-13pt

\noindent
{\sl same or opposite orientation.}

\medskip

{\bf Proof.} The process of finding such a decomposition of $\partial D$ for a graph-surface with one triangular face in a sphere is illustrated in the picture below. The resulting boundary polygon has $15$ edges.

The boundary $\partial D$ does not intersect the interior of any triangle $t\in T$. Let $\partial D\cap e\ni p$ be an inner point on an edge $e\in E$. The edge $e$ divides a small regular open neighbourhood
$U\ni p$ of $p$ into three disjoint parts $U=U'\sqcup(U\cap e)\sqcup U''$ with open $\varnothing\ne U'\subset D$ and $U''$.
Hence, $U\cap e\subset\partial D$. Thus, the boundary $\partial D$ contains the interior of $e$.

\medskip

\noindent
\includegraphics[width=0.35\textwidth]{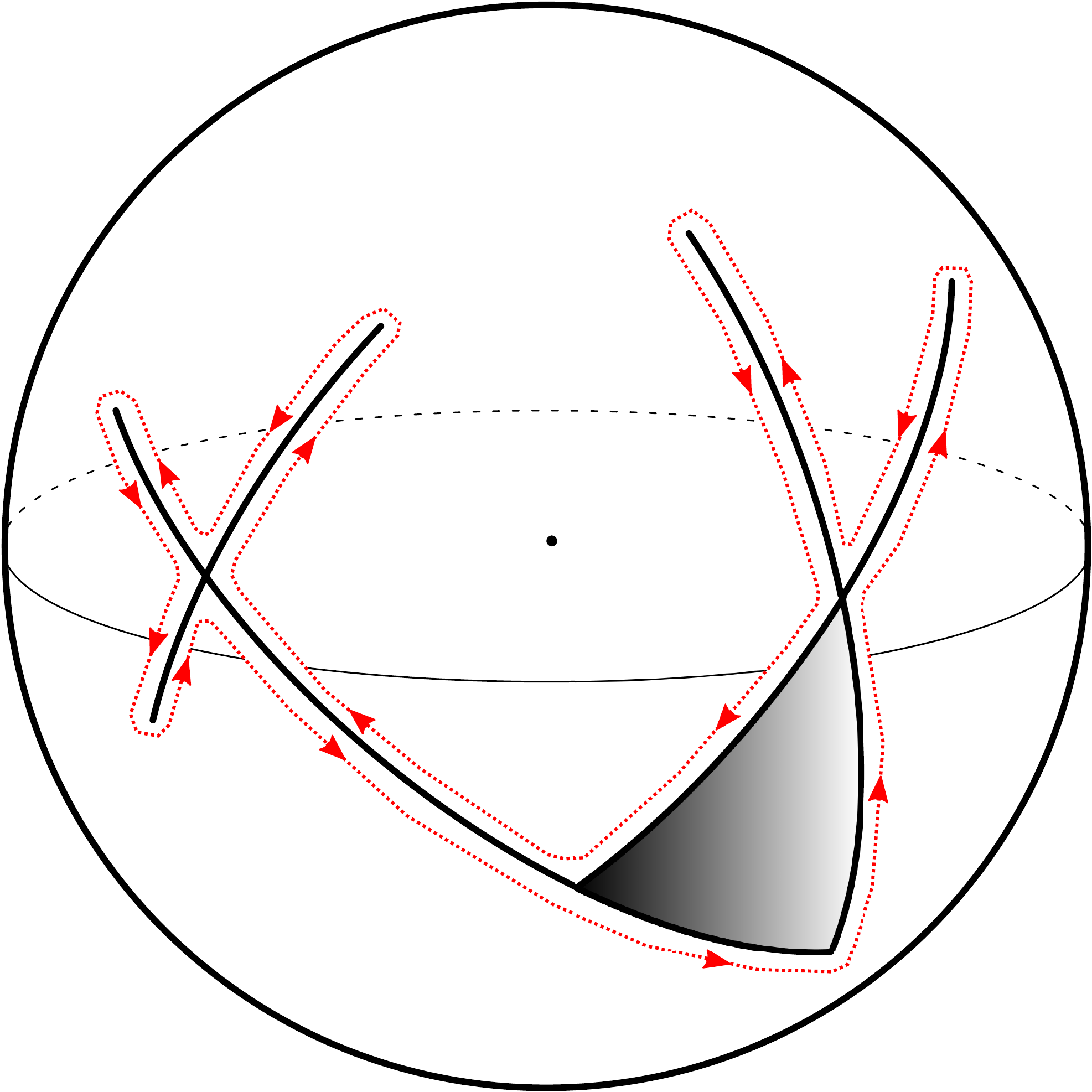}

\leftskip180pt

\vskip-160pt

Let $\partial D\ni v$ be a vertex, $v\in V$. For small $\varepsilon$, an open $\varepsilon$-neighbourhood $U\ni v$ of $v$ is divided into sectors by the
edges incident to $v$. Therefore, $D$ contains the interior of a sector, and, consequently, the interior of a suitable edge
$e\ni v$ intersects the boundary $\partial D$. Summarizing, we arrive at the first assertion.

For small $\varepsilon$, a closed $\varepsilon$-neighbourhood $Z\supset\partial D$ of the boundary $\partial D$ is a $2$-manifold with boundary. Let us denote this boundary by $B$. The
open disk $D$ contains exactly one component $C$ of $B$.
Taking, if necessary, a smaller $Z$, we induce a desired cyclic order from the circle $C$.~$_\blacksquare$

\leftskip0pt
\vskip20pt

\medskip

{\bf2.3.~Definition.} The edges from the list $g_1,\dots,g_k$ in Lemma 2.2 are called the {\it boundary edges\/} of the
graph-surface $S$.

\medskip

\noindent
\hskip260pt\includegraphics[width=0.23\textwidth]{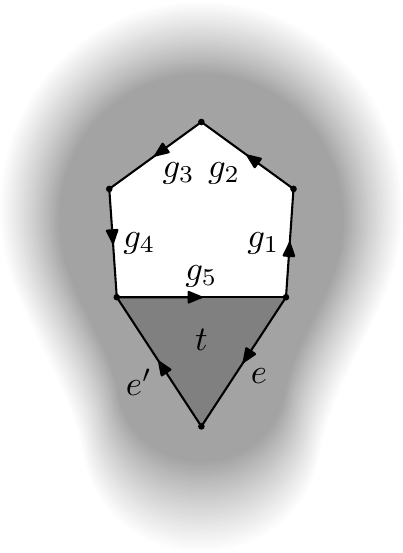}\includegraphics[width=0.23\textwidth]{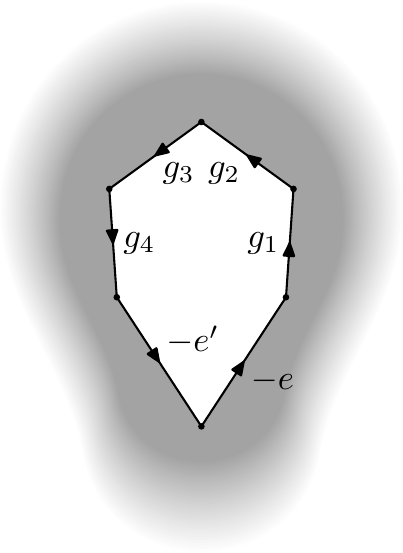}

\rightskip240pt

\vskip-150pt

{\bf2.4.~Collapse.} If a boundary edge $g_i$ is an edge of some triangle $t\in T$, then we can {\it collapse\/} the triangle
$t$, i.e., remove from the graph-surface $S$ the interior of the edge $g_i$ and the interior of the triangle $t$. The resulting
simplicial complex $S'\subset\hat S$ has the same vertices, $V'=V$, one less pair of oriented edges, $E'=E\setminus\{g_i,-g_i\}$,
and one less triangle, $T'=T\setminus t$. It is easy to see that $S'$ is again a graph-surface. The reader can readily
verify that everything works even if $t$ has two boundary edges.

\rightskip0pt

The decomposition $\partial D=g_1\cup\dots\cup g_i\cup\dots\cup g_k$ from Lemma 2.2 provides a similar decomposition
$\partial D'=g_1\cup\dots\cup(-e')\cup(-e)\cup\dots\cup g_k$ for the graph-surface $S'$, where the sequence $g_i$, $e$, $e'$ of
edges constitutes an oriented boundary $\partial t$ of $t$.

\medskip

{\bf2.5.~Remark.} If a graph-surface $S$ contains a triangular face, then it contains a face with a boundary edge (so that a collapse is
possible). Indeed, otherwise every boundary edge $g_i$ has the open disk $D$ ``on both sides'' of $g_i$, implying the inclusion
$S\subset\partial D$ and the absence of triangles in $S$.

\medskip

{\bf2.6.~Definition.} Let $P$ be a finite $1$-complex whose underlying space is a circle. So,~$P=U\cup F$, where $U$ is the set of
all vertices and $F=\{f_1,\dots,f_k\}$ the set of nondegenerate edges oriented and cyclically ordered with respect to the orientation of
the circle. Such a complex $P$ equipped with an edge length function $\ell:F\to]0,\infty[$ is called a {\it sample polygon.}

\medskip

{\bf2.7.~Remark.} Any closed surface $\hat S$ can be glued from a closed disk $\bar D$ as follows. We assume the
boundary $\partial D:=\bar D\smallsetminus D$, where $D \subset \bar D $ is the open disk in $\bar D$, to be the
underlying space of a finite $1$-complex, $\partial D=U\cup F$, where $U$ is the set of vertices and $F$ is the set of
nondegenerate edges taken with all orientations. Suppose that we are given an involution $\overline{\phantom f}:F\to F$ that commutes
with the orientation-flipping function $-:F\to F$, $f\mapsto-f$, that is $\overline{-f}=-\overline f$ for all $f\in F$. The
involution $\overline{\phantom f}$ indicates the pair $f,\overline f\in F$ of oriented edges to be glued.

Conversely, let $S\subset \hat S$ be a graph-surface without triangles. Then $\hat S$ is glued from a closed disk $\bar D$ and
$S$ coincides with the image of the boundary $\partial D$. For the proof, one can use the arguments explored in Remark 2.5 and in
the proof of Lemma 2.2.

Note that the list of edges $f_1,\dots,f_k$ in $F$ given by Lemma 2.2 provides such a gluing scheme $\overline{\phantom f}:F\to F$ for $\hat S$ when $S$ is a graph. Namely, $\overline{f_i}=f_j$ if the corresponding edges coincide up to sign $g_i=\pm g_j$. After the gluing, the boundary  $\partial D$ becomes the graph $S\subset \hat S$. In particular, if $\hat S$ is orientable, then the edges $f$ and $\overline f$ have different orientations with respect to the
orientation of the boundary $\partial D$ for any $f\in F$. In other words, when $S$ has no triangles and $\hat S$ is orientable,
each $g_i$ appears in a pair together with its opposite $-g_i$.

\medskip

{\bf2.8.~Kapovich-Millson-Klyachko's moduli space of polygons.} Pick a sample polygon $P$ and denote by $\mathbb E^P$ the set of all continuous maps
$P\to\mathbb E$ that are isometries on the edges. We call $\mathbb E^P$ the {\it space of polygons.} The {\it moduli space of
polygons\/} is the quotient $\mathbb E^P/\Isom^+\mathbb E$, where $\Isom^+\mathbb E$ is the group of all orientation-preserving isometries of $\mathbb E$. It carries the symplectic structure of Kapovich-Millson-Klyachko.

\medskip

{\bf2.9.~Explicit description of the space $\mathbb E^P/\mathbb E$.} Denoting by $\mathbb E\triangleleft\Isom^+\mathbb E$ the
normal subgroup of all translations, the quotient space $\mathbb E^P/\mathbb E$ is naturally identified with the set of all maps
$p:F\to\mathbb E$ such that
$$\langle p(f_i),p(f_i)\rangle=(\ell (f_i))^2\text{ for all }i,\hskip90pt\sum_ip(f_i)=0.\leqno{\mathbf{(2.9.1)}}$$
where $\langle-,-\rangle$ is the inner product in the euclidean linear space $\mathbb E$ and $\ell$ is the edge length function of a sample polygon.

Indeed, a continuous map $w:P\to\mathbb E$ that is isometric on edges provides an oriented segment $w(f_i)$ of length
$\ell (f_i)$ for every $i$. The associated vectors $p(f_i)\in\mathbb E$ in the linear euclidean space $\mathbb E$ obviously satisfy
the above identities.

Conversely, given a map $p:F\to\mathbb E$ that satisfies the identities, the map $w:P\to\mathbb E$ that sends the oriented edge
$f_j$ onto the oriented segment joining the points $\sum_{i=1}^{j-1}p(f_i),\sum_{i=1}^jp(f_i)\in\mathbb E$ obviously belongs to
$\mathbb E^P$.

\medskip

{\bf2.10.~Warning.} In spite of the description 2.9 being that simple, we actually deal with a real scheme over the real numbers
$\mathbb R$, i.e., with the space given by equations (2.9.1). At this stage we do not know yet whether the scheme is reduced. We prove that it is indeed reduced in Lemma 3.2.2.

\medskip

{\bf2.11.~Zariski tangent space to the scheme $\mathbb E^P/\mathbb E$.} Let $p:F\to\mathbb E$ be a point in $\mathbb E^P/\mathbb E$
(see~2.9). Taking derivatives of the equations (2.9.1) with respect to $p$, we obtain the equations for a tangent vector at $p$
to the scheme in question:
$${\T}_p(\mathbb E^P/\mathbb E)=\big\{t:F\to\mathbb E\mid\langle t(f_i),p(f_i)\rangle=0\text{ for all }i\text{ and
}\textstyle\sum_it(f_i)=0\big\}.$$

\medskip

{\bf2.12.~Definition.} A graph-surface equipped with an edge length function $\ell:E\to]0,\infty[$ such that $\ell (e)=\ell(-e)$
for all $e\in E$ and $\ell (e_1)+\ell (e_2)>\ell (e_3)$ for any triangle $t\in T$ whose boundary equals $\partial t=e_1+e_2+e_3$ (note
that we require three strict triangle inequalities for every triangle $t\in T$) is called a {\it sample polyhedron.} Let $S$ be a
sample polyhedron. Denote by $\mathbb E^S$ the set of all continuous maps $S\to\mathbb E$ that are isometries on the edges and
affine on the triangles. The~space $\mathbb E^S$ is called the {\it space of polyhedra.}

Given a sample polyhedron $S$, a graph-surface $S'$ obtained from $S$ by means of a collapse (see 2.4) is equipped with an edge
length function, the restriction $\ell'$ of the edge length function $\ell$ of $S$. So, $S'$ naturally becomes a sample
polyhedron.

\medskip

{\bf2.13.~Explicit description of the space $\mathbb E^S/\mathbb E$.} Let $u:S\to\mathbb E$ be a continuous map that is isometric
on the edges and affine on the triangles. For any oriented edge $e\in E$ (we remind the reader that, according to 2.1, $E$ is
the set of edges of $S$ taken with all orientations and equipped with the function $-:E\to E$, $e\mapsto-e$, that flips the
orientation), the image $u(e)$ is an oriented segment of length $\ell (e)$. Denote by $q(e)\in\mathbb E$ the associated vector in the
euclidean linear space $\mathbb E$. Then we get a map $q:E\to\mathbb E$ such that

\smallskip
\begin{equation}
\begin{minipage}{0.85\textwidth}

$\bullet$ $\langle q(e),q(e)\rangle=(\ell (e))^2$ for all $e\in E$,

$\bullet$ $q(-e)=-q(e)$ for all $e\in E$,

$\bullet$ $\sum_{i=1}^mq(e_i)=0$ for any closed path in $\cup E\subset S$ constituted by the edges $e_1,\dots,e_m\in E$.

\end{minipage}
\tag{\bf 2.13.1}
\end{equation}

\medskip

The next lemma simply states that a polyhedron can be reconstructed up to a parallel translation from the position of all its edges parallel translated to the origin of $\mathbb E$. Moreover, any set of vectors satisfying the conditions (2.13.1) gives rise to a unique up to a parallel translation polyhedron with directions and lengths of edges given by these vectors.

\medskip

{\bf2.14.~Lemma.} {\sl The quotient space\/ $\mathbb E^S/\mathbb E$ is naturally identified with the set of all maps\/
$q:E\to\mathbb E$ satisfying the above conditions.}

\medskip

{\bf Proof.} 
In the previous paragraph we constructed the map $q:\; E \to {\Bbb E}$
taking an edge $[a, b]= u(e)$ to $R(b)$, where
$R$ is the parallel translation taking $a$ to $0$.

On the other hand, given a map $q:\; E \to {\Bbb E}$
satisfying the above conditions (2.13.1), we define a map $u:\; S \to {\Bbb E}$
inductively. We start from an arbitrary edge $e_0\in S$,
and map it isometrically to the interval $[0, q(e_0)]$. The second edge $e_1$
is mapped to $[u(e_0), u(e_0)+q(e_1)]$, and so on, with
the edge $e_k$ taken to the interval
$\left[\sum_{i=0}^{k-1} q(e_i), \sum_{i=0}^{k} q(e_i)\right]$.

The only obstacle to defining the map $u$ on all the edges of $S$ in
this manner is the situation when the map $u$ is defined twice and differently on the same edge $e$. This means that we
have a closed path of edges inside $S$ that begins and ends with~$e$. It follows from the third condition from (2.13.1) that it can never happen and hence $u$ is
well-defined.

After having defined the map $u$ on $\cup E$, it is possible to obtain the unique extension of the map $u$ to $S$, affine on every triangle
$_\blacksquare$

\medskip

Every element $c\in\H_1(S,\mathbb Z)$ is representable by a closed path $\sum_{e\in E}h_ee$, where $h_e\in\mathbb Z$ for
all $e\in E$. Two homologous paths impose the same restrictions in the sense of 2.13, since they differ by a sum of boundaries of triangles. Thus it is enough to pick a set $H\subset\H_1(S,\mathbb Z)$ of generators of the abelian group $\H_1(S,\mathbb Z)$ and
rewrite the last condition in 2.13 as a couple of conditions

\smallskip

\noindent
$\bullet$ $qe_1+qe_2+qe_3=0$ for any triangle $t\in T$ with the boundary $\partial t=e_1+e_2+e_3$,

\noindent
$\bullet$ $\sum_{e\in E}h_eqe=0$ for any generator $\sum_{e\in E}h_ee\in H$.

\medskip

{\bf2.15.~Warning.} The space $\mathbb E^S/\mathbb E$ is described in 2.13 as a scheme. In the sequel we call the space
$\mathbb E^S/\mathbb E$ the {\it scheme of polyhedra.} Note that it may have nilpotents. 

\medskip

{\bf2.16.~Tangent space to the scheme $\mathbb E^S/\mathbb E$.} Let $q:E\to\mathbb E$ be a point in $\mathbb E^S/\mathbb E$ (see
2.13). Taking derivatives of the equations in 2.13, we obtain the equations for a tangent vector at $q$ to
the scheme $\mathbb E^S/\mathbb E$ :

\smallskip

The tangent space $\T_q(\mathbb E^S/\mathbb E)$ consists of all maps $s:E\to\mathbb E$ that satisfy the identities

\begin{equation}
\begin{minipage}{0.85\textwidth}
\noindent
$\bullet$ $\langle s(e),q(e)\rangle=0$ for all $e\in E$,

\noindent
$\bullet$ $s(-e)=-s(e)$ for all $e\in E$,

\noindent
$\bullet$ $s(e_1)+s(e_2)+s(e_3)=0$ for any triangle $t\in T$ with the boundary $\partial t=e_1+e_2+e_3$,

\noindent
$\bullet$ $\sum_{e\in E}h_es(e)=0$ for any representative $\sum_{e\in E}h_ee\in H$.
\end{minipage}
\tag{\bf 2.16.1}
\end{equation}

\medskip

{\bf2.17.~Boundary map.} Let $S$ be a sample polyhedron. We are going to  slightly rephrase the main statement of Lemma 2.2. For
any $1\leqslant i\leqslant k$, we introduce an oriented edge $f_i$ of length $\ell(g_i)$. Then we glue the edges $f_1,\dots,f_k$
into a sample polygon $P$ as in Definition 2.6. We~get a map $\delta:F\to E$, $\delta:f_i\mapsto g_i$, and a continuous map
$\bar\delta:P\to S$ induced by $\delta$. This means that $\bar\delta f_i=g_i$ for all $i$ and that the map $\bar\delta$ is
isometric on the edges. Clearly, $\bar\delta P=\partial D$. We call the map $\bar\delta:P\to S$ (or even the map $\delta:F\to E$)
the (combinatorial) {\it boundary map\/} of the sample polyhedron $S$.

The boundary map induces continuous maps $\mathbb E^S\to\mathbb E^P$ and $\mathbb E^S/\mathbb E\to\mathbb E^P/\mathbb E$ and
a linear map $\T_q(\mathbb E^S/\mathbb E)\to\T_{q\circ\delta}(\mathbb E^P/\mathbb E)$ for any point $q:E\to\mathbb E$ of the space
$\mathbb E^S/\mathbb E$ (see~2.13).

In terms of 2.9 and 2.13, the map $\mathbb E^S/\mathbb E\to\mathbb E^P/\mathbb E$ is given by the rule
$$(q:E\to\mathbb E)\mapsto(q\circ\delta:F\to\mathbb E).$$

In terms of 2.11 and 2.16, the map $\T_q(\mathbb E^S/\mathbb E)\to\T_{q\circ\delta}(\mathbb E^P/\mathbb E)$ is defined by a
similar rule
$$(s:E\to\mathbb E)\mapsto(s\circ\delta:F\to\mathbb E).$$

\bigskip

\centerline{\bf3. Symplectic structure on the moduli space of polygons}

\medskip

In this section we explicitly describe and study a natural skew-symmetric form on the space $\mathbb E^P/\mathbb E$ of polygons. After taking quotient by the group $G:=\SO(3,\mathbb R)$ this form descends to the Kapovich-Millson-Klyachko symplectic structure on the moduli space of polygons \cite{KaM,Klyachko}.
We concentrate on the linear-algebraic aspect of symplectic reduction. We refer to \cite{KaM} and standard references on symplectic geometry for the global questions.

\medskip

{\bf3.1.~Tangent space to a $G$-orbit.} Denote by $\so3$ the Lie algebra of the Lie group $G:=SO(3)$,
$$\so3:=\big\{a\in{\Hom}_\mathbb R(\mathbb E,\mathbb E)\mid\langle a(e),e'\rangle+\langle e,a(e')\rangle=0\text{ for all
}e,e'\in\mathbb E\big\}.$$
Let $p\in\mathbb E^P/\mathbb E$ be a point. By the definition  $p:F\to\mathbb E$ is a map such that
$\langle p(f_i),p(f_i)\rangle=(\ell (f_i))^2$ for all $i$ and $\sum_ip(f_i)=0$. The tangent space $\T_pGp$ to
the $G$-orbit of $p$ is the image of the Lie algebra $\so3$:
$${\T}_pGp=\{a\circ p:F\to\mathbb E\mid a\in \so3\},$$
where $a\circ p(f_i):=a(p(f_i))$ for all $i$.
Similarly, for a point $q\in\mathbb E^S/\mathbb E$, which is a map $q:E\to\mathbb E$ subject to the conditions listed
in 2.13, we get
$${\T}_qGq=\{a\circ q:E\to\mathbb E\mid a\in \so3\}.$$

\medskip

{\bf3.2.~Singular points of the scheme of polygons.} We discard from our consideration the cases where the length of an edge is
greater or equal than the sum of lengths of the other edges; in this case the moduli space of polygons has at most one point,
$|\mathbb E^P/\Isom^+\mathbb E|\leqslant1$. 

Taking into account that $t(f_1)=-\sum_{i=2}^kt(f_i)$, the remaining linear equations for a tangent vector $t:F\to\mathbb E$, that is,
for a sequence of vectors $t(f_i)\in\mathbb E$, $2\leqslant i\leqslant k$, take the form
$$\langle t(f_i),p(f_i)\rangle=0\text{ for all }2\leqslant i\leqslant k,\qquad\Big\langle\sum_{i=2}^kt(f_i),p(f_1)\Big\rangle=0.\leqno{(\mathbf{3.2.1})}$$
If these equations are linearly independent, then $\dim\T_p(\mathbb E^P/\mathbb E)=2k-3$, and the point $p$ is smooth. The first
$k-1$ linear equations are linearly independent because $p(f_i)\ne0$ for all $2\leqslant i\leqslant k$.

We claim that the last equation is not a linear combination of the first $k-1$ equations if the vectors $p(f_i)$,
$1\leqslant i\leqslant k$, are not collinear. Indeed, we can interpret the right-hand side of any equation in question as an element
in $\mathbb E\otimes_\mathbb R\mathbb E$, interpreting the $x_i:=t(f_i)$, $2\leqslant i\leqslant k$, as variables varying in
$\mathbb E$. Thus, the left-hand sides of equations correspond to $x_i\otimes p(f_i)$, $2\leqslant i\leqslant k$, and $\sum_{i=2}^kx_i\otimes p(f_1)$. Now
the claim that that the last equation from (3.2.1) follows from the rest is obvious.

We have proved the following
\medskip

{\bf3.2.1.~Lemma.} {\sl A point\/ $p\in\mathbb E^P/\mathbb E$, $p:F\to\mathbb E$, is singular if and only if all the vectors\/ $pf_i$ are
collinear.}~
$_\blacksquare$

\medskip

{\bf3.2.2.~Lemma.} {\sl If\/ $\mathbb E^P/\mathbb E\ni p$ is a singular point, then\/
$\dim\T_p(\mathbb E^P/\mathbb E)=1+\dim(\mathbb E^P/\mathbb E)$. The scheme\/ $\mathbb E^P/\mathbb E$ of polygons is reduced.}
\medskip

{\bf Proof.} The first assertion follows from the proof of Lemma 3.2.1: if all the vectors $pf_i$ are collinear then the last equation of the system of equations (3.2.1) follows from the first $k-1$ equations.

The second assertion follows from the fact that the scheme of polygons is a complete intersection. Indeed, the smooth locus has dimension $2k-3$, hence by \cite[Corollary 2.8.9]{BCR} this is the Krull dimension of $\mathbb E^P/\mathbb E$. On the other hand, the number of equations defining $\mathbb E^P/\mathbb E$ is
$k+3$. Thus it is a complete intersection that is generically smooth and hence generically reduced. Complete intersections
are Cohen-Macaulay, hence by \cite[Exercise 18.9]{Eisenbud} $\mathbb E^P/\mathbb E$ is reduced.
$_\blacksquare$

\medskip

{\bf3.3.~A skew-symmetric form on the scheme of polygons.} As in 2.9 and 2.11, the~tangent space to the product
$$\Pi:=\big\{\hat p:F\to\mathbb E\mid\langle\hat p(f_j),\hat p(f_j)\rangle=(\ell (f_j))^2\text{ for all }j\big\}$$
of $k$ $2$-spheres has the form
$${\T}_{\hat p}\Pi=\{\hat t:F\to\mathbb E\mid\langle\hat t(f_j),\hat p(f_j)\rangle=0\text{ for all }j\}.$$
As every $2$-sphere is endowed with its genuine symplectic form, the weighted sum 
$$\omega_{\hat p}(\hat t,\hat{t'}):=\sum_{j=1}^k\frac{\hat t(f_j)\wedge\hat{t'}(f_j)\wedge\hat p(f_j)}{(\ell(
f_j))^2\nu}\leqno{\mathbf{(3.3.1)}}$$
of these forms is a symplectic form on the product $\Pi$, where $\nu$ stands for the volume form on $\mathbb E$ and
$\hat t,\hat{t'}\in\T_{\hat p}\Pi$. This definition coincides with \cite[formula (b) below Remark 3.2, p.~491]{KaM}.

Since the scheme of polygons is a subscheme in $\Pi$, we get a skew symmetric form $\omega$ on the scheme of polygons given by
the same formula
$$\omega_p(t,t'):=\sum_{j=1}^k\frac{t(f_j)\wedge t'(f_j)\wedge pf_j}{(\ell (f_j))^2\nu}.\leqno{\mathbf{(3.3.2)}}$$

\medskip

{\bf3.3.3.~Remark.} Let $C$ be a finite-dimensional linear space equipped with a nondegenerate skew-symmetric form $\omega$, and
let $C\geqslant B$ be a linear subspace. Then $\dim\ker(\omega|_B)\leqslant\dim C-\dim B$. Indeed, otherwise the linear subspace
$\big\{b\in\ker(\omega|_B)\mid\omega(b,B')=0\big\}\subset\ker\omega$ would not be null, where $B'$ is a subspace complementary to
$B$, $B\oplus B'=C$, $\dim B'=\dim C-\dim B$.

\medskip

The next lemma is essentially a linear algebraic part of symplectic reduction applied to our set-up in an ad hoc manner. Note that the statement holds at regular and singular points.

\medskip

{\bf3.3.4.~Lemma.} {\sl The tangent space\/ $\T_pGp$ to the\/ $G$-orbit of any point\/ $p\in\mathbb E^P/\mathbb E$ coincides with
the kernel of the form\/ $\omega_p$ on $\T_p(\mathbb E^P/\mathbb E)$.}

\medskip

{\bf Proof.} First, we show that $\omega_p(t,a(p))=0$ for all $t\in\T_p(\mathbb E^P/\mathbb E)$ and $a\in \so3$, i.e., that
$\T_pGp\subset\ker\omega_p$.

Choosing a suitable orthonormal basis in $\mathbb E$, we may assume that
$a=\left[\begin{smallmatrix}0&1&0\\-1&0&0\\0&0&0\end{smallmatrix}\right]$. In this basis
$t(f_j)=\left[\begin{smallmatrix}t_{j1}\\t_{j2}\\t_{j3}\end{smallmatrix}\right]$ and
$p(f_j)=\left[\begin{smallmatrix}p_{j1}\\p_{j2}\\p_{j3}\end{smallmatrix}\right]$. Hence, $\sum\limits_{i=1}^3t_{ji}p_{ji}=0$ and
$\sum\limits_{i=1}^3p_{ji}^2=(\ell (f_j))^2$ for all $1\leqslant j\leqslant k$. Also,
$\sum\limits_{j=1}^kt_{ji}=\sum\limits_{j=1}^kp_{ji}=0$ for all $1\leqslant i\leqslant3$. It follows from
$\sum\limits_{i=1}^3t_{ji}p_{ji}=0$ and $\sum\limits_{i=1}^3p_{ji}^2=(\ell f_j)^2$ that
$$\frac{t(f_j)\wedge a\circ p(f_j)\wedge p(f_j)}\nu=\frac{\left[\begin{matrix}t_{j1}\\t_{j2}\\t_{j3}\end{matrix}\right]\wedge
\left[\begin{matrix}p_{j2}\\-p_{j1}\\0\end{matrix}\right]\wedge\left[\begin{matrix}p_{j1}\\p_{j2}\\p_{j3}\end{matrix}\right]}
\nu=\det\left[\begin{matrix}t_{j1}&p_{j2}&p_{j1}\\t_{j2}&-p_{j1}&p_{j2}\\t_{j3}&0&p_{j3}\end{matrix}\right]=$$
$$=-t_{j1}p_{j1}p_{j3}-t_{j2}p_{j2}p_{j3}+t_{j3}p_{j2}^2+t_{j3}p_{j1}^2=t_{j3}p_{j3}^2+t_{j3}p_{j2}^2+t_{j3}p_{j1}^2=
t_{j3}(\ell f_j)^2.$$
It remains to use the equality $\sum\limits_{j=1}^kt_{j3}=0$.

Next, we take $C:=\T_p\Pi$ and $B:=\T_p(\mathbb E^P/\mathbb E)$  as in Remark 3.3.3. From $\dim\T_p\Pi=2k$,
$\dim(\mathbb E^P/\mathbb E)=2k-3$, Lemma 3.2.2, and the first part of the proof of Lemma~3.3.4, we~conclude that
$\dim C-\dim B=3$ for a smooth point $p$, that $\dim C-\dim B=2$ for a singular point $p$, and that
$\dim\T_pGp\leqslant\dim\ker\omega_p$. Using the description 3.1 of the tangent space to a $G$-orbit, it remains to observe that
$\dim\T_pGp=3$ for a smooth point $p$ and $\dim\T_pGp=2$ for a singular point $p$
$_\blacksquare$

\bigskip

\centerline{\bf4.~Proofs of main results}

\medskip

The following trivial lemma is a classical fact claiming that any infinitesimal deformation of a triangle extends to an infinitesimal rotation of $\mathbb E$. In other words, any Zariski tangent vector at a polygon $p$, consisting of a single triangular face, belongs to $\T_pGp$.

{\bf4.1.~Lemma {\rm(the rigidity of a triangle)}.} {\sl Let\/ $p_1,p_2,p_3\in\mathbb E$ be noncollinear and such that\/
$p_1+p_2+p_3=0$. Suppose that\/ $\langle t_j,p_j\rangle=0$ for all\/ $j$ and\/ $t_1+t_2+t_3=0$, where\/
$t_1,t_2,t_3\in\mathbb E$. Then there exists an element\/ $a\in \so3$ such that\/ $t_j=a(p_j)$ for all\/ $j$.}

\medskip

{\bf Proof.} We take a linear map $a\in\Hom_\mathbb R(\mathbb E,\mathbb E)$ such that $a:p_1\mapsto t_1$ and $a:p_2\mapsto t_2$.
It~follows from $p_1+p_2+p_3=0$ and $t_1+t_2+t_3=0$ that $a:p_3\mapsto t_3$. Since $\langle t_j,p_j\rangle=0$, we obtain
$\langle a(p_j),p_j\rangle=0$ for all $j$. For pairwise distinct indices $j,k,l$, we have
$$\langle a(p_j),p_k\rangle+\langle p_j,a(p_k)\rangle=\langle t_j,p_k\rangle+\langle p_j,t_k\rangle=-\langle
t_k+t_l,p_k\rangle-\langle p_j,t_j+t_l\rangle=$$
$$=-\langle t_l,p_k\rangle-\langle p_j,t_l\rangle=-\langle t_l,p_k+p_j\rangle=\langle t_l,p_l\rangle=0.$$
We pick $p\ne0$ such that $\langle p,p_j\rangle=0$ for all $j$ and define a unique $a(p)\in\Span(p_1,p_2,p_3)$ subject to
$\langle a(p),p_j\rangle=-\langle p,t_j\rangle$ for all $j$ (any two of these three equalities imply the third one). Then
$\langle a(p),p\rangle=0$ and $\langle a(p),p_j\rangle+\langle p,a(p_j)\rangle=0$ for all $j$. Hence, $a\in \so3$.
$_\blacksquare$

\medskip

\begin{wrapfigure}{l}{0.5\textwidth}
  \begin{center}
  \begin{tikzcd}
\T_{q}\mathbb E^{S}/\mathbb E \arrow{r}{\text{restriction}}[swap]{\text {to }S'\subset S} \arrow[swap]{d}{d\delta} & \T_{q'}\mathbb E^{S'}/\mathbb E \arrow{d}{d\delta'} \\%
\T_{\delta(q)}\mathbb E^{P}/\mathbb E,\omega & \T_{\delta'(q)}\mathbb E^{P'}/\mathbb E, \omega' \\
\arrow[hookrightarrow]{u}{\text{inclusion}} d\delta(\T_{q}\mathbb E^{S}/\mathbb E ) \arrow[hookrightarrow,dashed]{r}{i} &
\arrow[hookrightarrow]{u}{\text{inclusion}} d\delta'(\T_{q'}\mathbb E^{S'}/\mathbb E )
\end{tikzcd}
  \end{center}
\caption*{Collapsing and the boundary map}
\end{wrapfigure}

The next proposition will be the core of our argument. Its meaning is the following. Suppose that one has a pair of Zariski tangent vectors $s_1,s_2\in \T_q(\mathbb E^S/\mathbb E)$ at a point $q\in \mathbb E^S/\mathbb E$ in the space of polyhedra associated to a sample polyhedron $S$. There is a point $q'\in \mathbb E^{S'}/\mathbb E$ in the space of polyhedra associated to a sample polyhedron $S'$, obtained from $S$ by the collapse of a face --- namely, one just removes the face from the polyhedron corresponding to $q$. This inclusion gives a pair of tangent vectors $s_1',s_2'\in \T_{q'}\mathbb E^{S'}/\mathbb E$ at $q'$ by restriction to the remaining edges (this is the first part of Proposition 4.2).
The second part of Proposition 4.2 states that if $\omega'(s_1',s_2')=0$, then $\omega(s_1,s_2)=0$.

We remind (2.17), that each sample polyhedron $S$ has {\it its own} associated space of polyhedra $\mathbb E^S/\mathbb E$, the boundary polygon $P$ with the corresponding space of polygons $\mathbb E^P/\mathbb E$ equiped with the skew-symmetric form $\omega$ (defined in 3.3), and the map $\delta$. 
The relation between the spaces involved is shown in the diagram on the left.

\medskip

{\bf4.2.~Proposition.} {\sl Let\/ $S$ and\/ $S'$ be sample polyhedra such that\/ $S'$ is made from\/ $S$ by means of the collapse
of a triangle, and let\/ $ q \in \mathbb E^S/\mathbb E$ be a point, $q:E\to\mathbb E$. Denote by\/ $q':E'\to\mathbb E$ the
restriction of\/ $q$ to\/ $E'\subset E$ and by\/ $\delta:F\to E$ and\/ $\delta':F'\to E'$ the corresponding boundary maps. Then

\noindent (I) $q'\in\mathbb E^{S'}/\mathbb E$

\noindent (II) Suppose that\/ $\omega_{q'\circ\delta'}(s'_1\circ\delta',s'_2\circ\delta')=0$ for all\/
$s'_1,s'_2\in\T_{q'}(\mathbb E^{S'}/\mathbb E)$. Then\/ 
 $\omega_{q\circ\delta}(s_1\circ\delta,s_2\circ\delta)=0$ for all\/
$s_1,s_2\in\T_q(\mathbb E^S/\mathbb E)$.}

\medskip

{\bf Proof.} 

{\bf Step 1.}
In order to show that $q'\in\mathbb E^{S'}/\mathbb E$, it suffices to verify the following four conditions equivalent to conditions (2.13.1):

\smallskip

\noindent
$(i)$ $\langle q'(e'),q'(e')\rangle=(\ell'(e'))^2$ for all $e'\in E'$,

\noindent
$(ii)$ $q'(-e')=-q'(e')$ for all $e'\in E'$,

\noindent
$(iii)$ $q'(e')_1+q'(e')_2+q'(e')_3=0$ for any triangle $t'\in T'$ with the boundary $\partial t'=e'_1+e'_2+e'_3$,

\noindent
$(iv)$ $\sum_{e'\in E'}h_{e'}q'(e')=0$ for any generator $\sum_{e'\in E'}h_{e'}e'\in H'$.

\smallskip

\noindent
for some choice of a set $H'\subset\H_1(S',\mathbb Z)$ of generators.

Since $\ell',q'$ are restrictions of $\ell,q$ and $T'\subset T$, the first three are immediate. For the fourth condition, we
take any set of paths in $S'$ generating $H_1(S',\mathbb{Z})\cong H_1(S,\mathbb{Z})$. Thus $(I)$ is proved.

{\bf Step 2.}
Next we are going to show that the restriction of any $s\in \T_q(\mathbb E^S/\mathbb E)$ to the set of edges of $S'$ gives a tangent vector to $\mathbb E^{S'}/\mathbb E$.

Let $s\in\T_q(\mathbb E^S/\mathbb E)$ or, equivalently, a map $s:E\to\mathbb E$ satisfying the identities (2.16.1). In order to prove that the restriction $s':E'\to\mathbb E$ of $s$ to $E'\subset E$ belongs to the tangent space
$\T_{q'}(\mathbb E^{S'}/\mathbb E)$, it suffices to verify the identities (2.16.1):

\smallskip

\noindent
$(i)$ $\langle s'(e'),q'(e')\rangle=0$ for all $e'\in E'$,

\noindent
$(ii)$ $s'(-e')=-s'(e')$ for all $e'\in E'$,

\noindent
$(iii)$ $s'(e')_1+s'(e')_2+s'(e')_3=0$ for any triangle $t'\in T'$ with the boundary $\partial t'=e'_1+e'_2+e'_3$,

\noindent
$(iv)$ $\sum_{e'\in E'}h_{e'}s'(e')=0$ for any generator $\sum_{e'\in E'}h_{e'}e'\in H'$
\smallskip

The first three are immediate because $T'\subset T$, the maps $s',q'$ are restrictions of the maps $s,q$, which in turn
satisfy the identities (2.16.1). The last identity holds automatically since the inclusion $S'\subset S$ induces an isomorphism $H'\cong H$.

{\bf Step 3.}
Now let $s_1,s_2\in\T_q(\mathbb E^S/\mathbb E)$. By Paragraph 3.1 and Lemma 3.3.4, to prove that $\omega_{q\circ\delta}(s_1\circ\delta,s_2\circ\delta)=0$ it suffices to show that
$$\omega_{q\circ\delta}\big((s_1-a_1(q))\circ\delta,(s_2-a_2(q))\circ\delta\big)=0$$
for some $a_1,a_2\in \so3$ because the differential map
$d\delta=(\cdot\circ \delta): \T_q(\mathbb E^S/\mathbb E)\to\T_{q\circ\delta}(\mathbb E^P/\mathbb E)$ sends the tangent space to a $SO(3)$-orbit to the tangent
space to a $SO(3)$-orbit and adding a vector tangent to an orbit does not change the value of the form (Lemma 3.3.4).

Recall that the sample polyhedron $S'$ is made from the sample polyhedron $S$ by the collapse of a triangle $t\in T$. As in paragraph 2.4,
we denote by $g_i$, $e$, $e'$ the sequence of edges constituting an oriented boundary $\partial t$ of $t$, where $g_i$ is a
boundary edge which is removed.

Restricting the map $s_j:E\to\mathbb E$ to the boundary $\partial t$ of the triangle $t$, i.e., to
$E_0:=\{-g_i,g_i,-e,e,-e',e'\}\subset E$, and taking $p_1:=q(g_i)$, $p_2:=q(e)$, $p_3:=q(e')$, $t_1:=s_j(g_i)$, $t_2:=s_j(e)$, $t_3:=s_j(e')$
in Rigidity of Triangle Lemma 4.1 (this lemma is applicable since $q$ satisfies (2.13.1) and  $s_j$ satisfies (2.16.1)), we find an element $a_j\in \so3$ such that $a_j\circ q(g_i)=s_j(g_i)$, $a_j\circ q(e)=s_je$, $a_j\circ q(e')=s_j(e')$ for $j=1,2$. This means that
$s_j-a_j(q)$ is null on $E_0$ for $j=1,2$. In other words, we may assume that 
\[
s_1 \text{ and } s_2 \text{ are null on } E_0
\leqno{\mathbf{(4.2.1)}}
\]

{\bf Step 4.} By (4.2.1) we have tangent vectors $s_1,s_2\in\T_q(\mathbb E^S/\mathbb E)$,
$s_1,s_2:E\to\mathbb E$, that vanish on $E_0:=\{-g_i,g_i,-e,e,-e',e'\}\subset E$. Their restrictions $s'_1,s'_2$ to $E'\subset E$
belong to the tangent space $\T_{q'}(\mathbb E^{S'}/\mathbb E)$ by Step 2. Now 
$\omega_{q'\circ\delta'}(s'_1\circ\delta',s'_2\circ\delta')=0$ by the assumptions of Proposition 4.2. The boundary edges of $S$
are listed in $\partial D=g_1\cup\dots\cup g_i\cup\dots\cup g_k$ in their cyclic order (as in Lemma 2.2). The boundary edges of $S'$ are
listed in $\partial D'=g_1\cup\dots\cup(-e')\cup(-e)\cup\dots\cup g_k$ in their cyclic order (as in Paragraph 2.4). Take the sets of edges of the correspondng combinatorial boundary polygons of $S$ and $S'$:
$F=\{f_1,\dots,f_i,\dots,f_k\}$ and $F'=\{f_1,\dots,f'_i,f''_i,\dots,f_k\}$. Now we have $\delta f_j=\delta'f_j=g_j$ for any
$j\ne i$ with $\delta f_i=g_i$, $\delta'f'_i=-e'$, and $\delta'f''_i=-e$.

Expanding the definition (3.3.2) for $\omega_{q\circ\delta}(s_1\circ\delta,s_2\circ\delta)$ and
$\omega_{q'\circ\delta'}(s'_1\circ\delta',s'_2\circ\delta')$, we~can see that almost all summands coincide because
$\delta (f_j)=\delta'(f_j)$ for all $j\ne i$ and $q',s'_1,s'_2$ are the restrictions of $q,s_1,s_2$. On the other hand, the remaining
summands are
$$\frac{(s_1\circ\delta)(f_i)\wedge(s_2\circ\delta)(f_i)\wedge(q\circ\delta)(f_i)}{(\ell (f_i))^2\nu}$$
for $\omega_{q\circ\delta}(s_1\circ\delta,s_2\circ\delta)$ and
$$\frac{(s'_1\circ\delta')(f'_i)\wedge(s'_2\circ\delta')(f'_i)\wedge(q'\circ\delta')(f'_i)}{(\ell
(e'))^2\nu}+\frac{(s'_1\circ\delta')(f'')_i\wedge(s'_2\circ\delta')(f''_i)\wedge(q'\circ\delta')(f''_i)}{(\ell (e))^2\nu}$$
for $\omega_{q'\circ\delta'}(s'_1\circ\delta',s'_2\circ\delta')$. A closer look allows us to infer that all these summands
vanish since $s_1$ vanishes on $E_0$, thus finishing the proof of (II).
$_\blacksquare$

\medskip

{\bf4.3.~Theorem.} {\sl Let\/ $\hat S\supset S$ be a graph-surface in an orientable closed surface. Then
$$\delta:\mathbb E^S/{\Isom}^+\mathbb E\to\mathbb E^P/{\Isom}^+\mathbb E$$
is isotropic. That is, the pullback of the form\/ $\omega$ is null.}

\medskip

{\bf Proof.}  
By Remark 2.5, after a series of collapses one ends up with a graph-surface %
without triangles. For such graph-surfaces on an orientable surface $\hat S$, by Remark 2.7, %
each edge $g_i$ of the combinatorial boundary always appears together with its opposite $-g_i$ in the decomposition of the boundary $\partial D$ %
given by Lemma 2.2. This immediately implies that the skew-symmetric form 
$$\omega_p(s_1\circ\delta,s_2\circ\delta)=\sum_{j=1}^k\frac{(s_1 \circ\delta) (f_j)\wedge (s_2\circ \delta) (f_j)\wedge p(f_j)}{(\ell (f_j))^2\nu}$$
is zero for any $s_1,s_2\in\T_q(\mathbb E^S/\mathbb E)$. In other words, for graph-surfaces without triangles the boundary %
map $\delta: \mathbb E^S/\mathbb E\to\mathbb E^P/\mathbb E$ is isotropic.

On the other hand, by Proposition 4.2, if  $\omega_{q'\circ\delta'}(s'_1\circ\delta',s'_2\circ\delta')=0$ for all %
$s'_1,s'_2\in\T_{q'}(\mathbb E^{S'}/\mathbb E)$ then \/ $\omega_{q\circ\delta}(s_1\circ\delta,s_2\circ\delta)=0$ %
for all $s_1,s_2\in\T_q(\mathbb E^S/\mathbb E)$, where $S'$ is obtained by a collapse of $S$. Hence the map %
$\delta: \mathbb E^S/\mathbb E\to\mathbb E^P/\mathbb E$ is indeed isotropic for arbitrary graph-surfaces $S$. %
The property of being isotropic survives in the quotient of $\T_q(\mathbb E^S/\mathbb E)$ and $\T_{q\circ \delta}(\mathbb E^P/\mathbb E)$ %
by the subspace tangent to the orbit of $SO(3)$ through $q$. %
$_\blacksquare$ %

\medskip

\noindent
\includegraphics[width=0.24\textwidth]{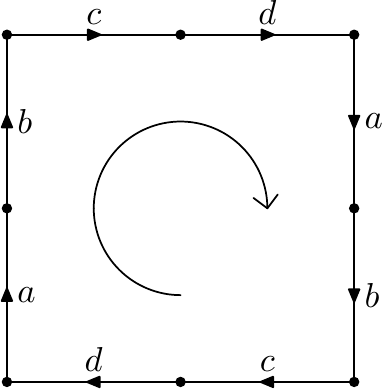}

\leftskip125pt

\vskip-114pt

{\bf4.4.~Remark.} If the surface $\hat S$ is not orientable, the statement of
the theorem 4.3 need not hold for a graph-surface. 
For example, consider
$\hat S=\mathbb {RP}^2$ obtained by identifying the sides of an 8-gon according to the scheme depicted on the picture on the left. The boundary of the 8-gon maps to the graph $S\subset \hat S$ isomorphic to the boundary of a square. Let us assign the lenght $1$ to all edges of $S$, that is $\ell(a)=\ell(b)=\ell(c)=\ell(d)=1$.
This is a sample polyhedron (without triangular faces) with the set of edges $\{\pm a,\pm b,\pm c,\pm d\}$ all of length $1$ living in the real  projective plane.

\leftskip0pt

\vskip15pt

Choosing an orientation on the complement disk and applying Lemma 2.2, as usual, we find a boundary polygon  $P=abcdabcd$. Note that each edge appears twice {\it with the same orientation} --- this is an artefact of non-orientability of $\hat S$.

Consider an embedding $q$ of $S$ to $\mathbb R^3$ and two Zariski tangent vectors $s_1,s_2\in \T_q(\mathbb E^S/\mathbb E)$
given by 

\[q= \begin{cases} 
      a\mapsto(0,-1,0)  \\ b\mapsto(1,0,0) \\
      c\mapsto(0,1,0) \\   d\mapsto(-1,0,0)
   \end{cases}
   s_1=\begin{cases} 
      a\mapsto(-1,0,0) \\
      b\mapsto(0,0,0) \\
      c\mapsto(1,0,0) \\
      d\mapsto(0,0,0)
   \end{cases}
  s_2=\begin{cases} 
      a\mapsto(0,0,1) \\
      b\mapsto(0,0,-1) \\
      c\mapsto(0,0,1) \\
      d\mapsto(0,0,-1)
   \end{cases}
\]

\medskip

One can easily check that $q$ indeed satisfies (2.13.1) and $s_1,s_2$ satisfy (2.16.1).

\smallskip

Now we can compute the value $\omega(s_1\circ\delta,s_2\circ\delta)$ at $q$:
$$\omega(s_1\circ\delta,s_2\circ\delta)=\sum_{e \in P} \langle s_1(e)\times s_2(e), q(e)\rangle=1+0+1+0+1+0+1+0=4$$
That is, the pull-back of $\delta$ is not null. This example cannot be made into a genuine triangulated surface by completing the edges to triangles. Hovever, we propose the following 

\medskip

{\bf4.5.~Conjecture.}  There exists a non-orientable graph-surface that is  topologically a surface with boundary, such that $\delta$ is not isotropic. 

\bigskip

\centerline{\bf5.~Co-isotropy of $\delta$}

\bigskip

\noindent
\includegraphics[width=0.24\textwidth]{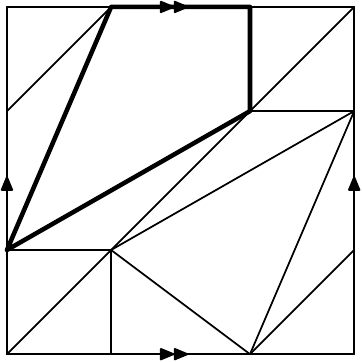}

\leftskip125pt

\vskip-120pt
This section is devoted to finding the {\it lower bound} on the dimension of the image of $\delta$. 

\medskip

{\bf5.1.~Remark.}
If the Euler characteristic of a surface is less than one, then the image of $\delta$ is generally strictly less than half the dimension of $\mathbb E^P/{\Isom}^+\mathbb E$. For example, the surface of the combinatorial type depicted in the figure on the left is generally infinitesimally rigid (see \cite{Korshunov} for a script that computes it). The boundary is shown in bold.

\leftskip0pt

\vskip25pt

{\bf5.2.~Lemma.} {\sl Let $S\subset \hat S$ be a graph-surface. Then $3|T|=2|E|-|F|$, where $|T|,|E|,|F|$ is a number of triangles, edges, and the number of edges in the boundary respectively.}

\medskip

{\bf Proof.} Any triangle has exactly $3$ adjacent edges. Any non-boundary edge belongs to two triangles. On the other hand, a boundary edge can be adjacent to one or no triangles. In the latter case there are $2$ elements of the boundary F (in the sense of Definition 2.3) corresponding to it, in the former --- one element of the boundary $F$ and one triangle. Hence by a simple inclusion-exclusion argument we obtain $3|T|+|F|=2|E|$. Indeed, each triangle gives $3$ edges and each element of the boundary $F$ gives one edge, and after summation one has each edge counted exactly twice.
$_\blacksquare$

\medskip

{\bf5.3.~Lemma.} {\sl Let $\mathbb E^S/\mathbb E$ be a configuration space of a sample polyhedron as in Definition 2.12. The dimension of $\mathbb E^S/\mathbb E$ is at least $|F|-\mathrm{dim}\H_1(S,\mathbb Z)$.}

\medskip

{\bf Proof.}
Note that by Lemma 2.14 we have that $\mathbb E^S/\mathbb E$ is given by $|E|+3|T|+\mathrm{dim}\,\H_1(S,\mathbb Z)$ equations in $3|E|$ variables. Thus we obtain a lower bound on the dimension:¨ $3|E|-|E|-3|T|-\mathrm{dim}\,\H_1(S,\mathbb Z)=2|E|-3|T|-\mathrm{dim}\,\H_1(S,\mathbb Z)$. Now, by Lemma 5.2, one has  $\dim \mathbb E^S/\mathbb E \ge 3|T|+|F|-3|T|-\mathrm{dim}\,\H_1(S,\mathbb Z)=|F|-\mathrm{dim}\,\H_1(S,\mathbb Z)$.
$_\blacksquare$

\medskip

From now on we suppose that our graph-surface $S$ is homeomorphic to a disk. To obtain a lower bound on the {\it image} of $\delta$ we need to establish an {\it upper bound} on the kernel of $\delta$.  
It is known \cite{Gluck}, that the set of infinitesimally rigid polyhedra of a given combinatorial type, homeomorhic to a sphere, is a Zariski open dense subset of $\mathbb R^{3|V|}$:

\medskip
{\bf5.4.~Proposition (Gluck's theorem).} {\sl Let $\hat S$ be a triangulated surface homeomorphic to a sphere. We will refer to functions from the set of vertices of $\hat S$ to $\mathbb R^3$ as polyhedra in $\mathbb R^3$ of a fixed combinatorial type $\hat S$ (possibly with degenerate edges). The set of all polyhedra of type $\hat S$ is naturally identified with $\mathbb R^{3|V|}$. Each polyhedron $q$ gives rise to a sample polyhedron in the sense of Definition 2.12 and the corresponding space of polyhedra $\mathbb E^{\hat S}$. Let $D\subset \mathbb R^{3|V|}$ be the set of all polyhedra $q$ for which the corresponding $\dim\,\T_q\mathbb E^S\ge 6$. Gluck's theorem states that $D$ comprises a proper algebraic subset of $\mathbb R^{3|V|}$. In other words, a generic polyhedron homeomorphic to a sphere is infinitesimally rigid (i.e. all its infinitesimal deformations are trivial).}

\medskip
{\bf A sketch of the Gluck's proof.} Consider the map 
$\hat\rho: \mathbb R^{3|V|} \to \mathbb R^{|E|}$ that sends a polyhedron to the $|E|$-tuple of the squares of lenghts of its edges:
$$\hat\rho: (\dots x_i,y_i,z_i,\dots x_j,y_j,z_j,\dots)\mapsto (\dots (x_i-x_j)^2+(y_i-y_j)^2+(z_i-z_j)^2,\dots)$$
for all adjacent $i$ and $j$. A fiber of $\hat\rho$ is exactly $\mathbb E^{\hat S}$ for a triangulated surface $\hat S$ with a {\it fixed metric} (a sample polyhedron in terms of Definition 2.12). The affine subspace $\mathrm{Ker}\, \hat\rho_*\subset \mathbb{R}^{3|V|}$ can be identified with the Zariski tangent space $\T_p\mathbb E^S$ at the corresponding polyhedron $p$. The dimension of a fiber $\hat\rho^{-1}(m)$ is at least $6$ (corresponding to the trivial deformations). The condition that the dimension of a fiber at a point is greater than $6$ is algebraic. Indeed, it is equivalent to checking that the kernel of the differential has dimension greater than $6$. This, in turn, is equivalent to the vanishing of all minors of rank greater than $3|V|-6$, which is an algebraic condition. Thus $D$ is an algebraic subset of $\mathbb{R}^{3|V|}$. 

It remains to prove that $D$ is a proper subset. By a theorem of Steinitz \cite[Section 23]{Lyusternik} any triangulated surface homeomorphic to the sphere admits a realization $q$ as a {\it convex} polyhedron in $\mathbb R^3$. By a theorem of Cauchy \cite[Section 20]{Lyusternik} it is infinitesimally rigid and hence $q$ does not lie in $D$. $_\blacksquare$

\medskip

{\bf5.5.~Definition.} We will say that a polyhedron (with boundary) $q\in \mathbb R^{3|V|}$ is {\it boundary rigid} if for any tangent vector $s\in  \T_q\mathbb E^S/\mathbb E$ the condition $d\delta(s)=0$  implies $s=0$. In plain language, after fixing the boundary, a polygon becomes infinitesimally rigid.

\bigskip
\noindent
\includegraphics[width=0.24\textwidth]{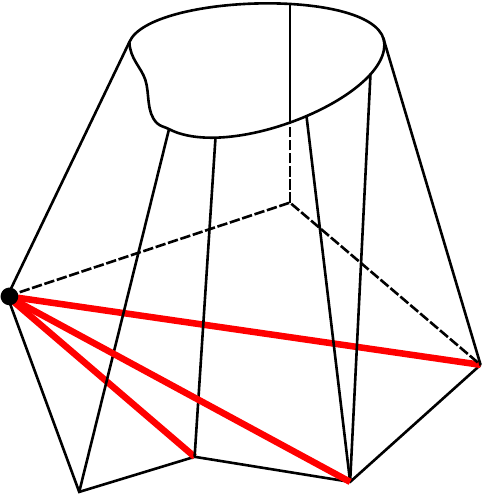}

\leftskip125pt

\vskip-123pt

{\bf5.6.~Remark.} Let $\hat S$ be a closed surface obtained from a surface with boundary $S$ by choosing a boundary vertex and adding all $|F|-3$ segments connecting the vertex with non-adjacent boundary vertices as shown in the picture on the left. The surface $\hat S$ is homeomorphic to the sphere. Note, that $S$ and $\hat S$ have the same set of vertices $V$. There is a natural identification of spaces of all polyhedra of combinatorial types $S$ and $\hat S$ with the space of all maps $V\to \mathbb R^3$ (i.e. with $\mathbb R^{3|V|}$).
 If a polyhedron with boundary $q$ of type $S$ admits an infinitesimal  deformation $s$ such that $d\delta(s)=0$, then $s$ also gives rise to an infinitesimal deformation $\hat s$ of a polyhedron  $\hat q$ of type $\hat S$: define $\hat s:=s$ on old edges and $\hat s:=0$ on added edges. The Conditions (2.16.1) obviously 
continue to hold. Thus we arrived at the following

\leftskip0pt

\vskip10pt

{\bf5.7.~Proposition.} {\sl If $\hat S$ is infintesimally rigid, then $S$ is boundary rigid. }$_\blacksquare$

\medskip

We denote by $\Delta$ the subset of $\mathbb R^{3|V|}$ corresponding to polyhedra of type $\hat S$ that are not infinitesimally rigid. In particular, $\Delta$ contaians all polyhedra of type $S$ that are not boundary rigid.

\medskip

We summarize the relationships between different configuration spaces and maps defined so far in the following diagram:

\medskip

\begin{center}
\begin{tikzcd}
\Delta\arrow[hookrightarrow]{r}{\subset} & \mathbb R^{3|V|} \arrow{r}{\rho} \arrow{rd}{\hat \rho} \arrow[swap]{d}{\delta} & \mathbb R^{|E|} \\
& \mathbb R^{3|F|} & \mathbb R^{|\hat E|}
\end{tikzcd}
\end{center}

\medskip

Here the maps $\rho$ and $\hat \rho$ are square-lenght maps as in the proof of the Gluck's theorem: a polyhedron of combinatorial type $S$ (respectively $\hat S$) is maped to the function $E\ni e \mapsto (\ell(e))^2$ (respectively $\hat E\ni \hat e \mapsto (\ell(\hat e))^2$). Thus a fiber of $\rho$ (respectively $\hat\rho$) is the space $\mathbb E^{S}$ (respectively $\mathbb E^{\hat S}$) of Definition 2.12, that is, the set of all polyhedra of a given combinatorial type and a {\it fixed inner metric}. The map $\delta$ sends a polyhedron with boundary of combinatorial type $S$ to a tuple of vectors constituting its boundary. On a fiber of $\rho$ it restricts to the $\delta$ of Paragraph 2.17.

\medskip

{\bf5.8.~Remark.} Note that Gluck's theorem only guarantees that a generic polyhedron {\it among all} polyhedra of combinatorial type $S$ is infinitesimally rigid. When we fix an inner metric on $\hat S$, that is, chose a fiber of $\hat\rho$, the resulting configuration space of isometric realizations $\mathbb E^S$ could contain a connected component which is non-rigid. The Gluck's argument does not address directly the question of rigidity of {\it all} components of a {\it general fiber}. 

The following lemma ensures that {\it every} connected component of a general fiber of $\rho$ contains a boundary rigid polyhedron:

\medskip

{\bf5.9.~Lemma.} {\sl Let $ \rho: \mathbb R^{3|V|}\to \mathbb R^{| E|}$ be the square-length map for a polyhedral surface homeomorphic to a disk. Then there is a subset $R\subset \rho(\mathbb R^{3|V|})$ of full measure such that for any $g\in R$ the intersection $\rho^{-1}(g)\cap \Delta$ is nowhere dense in $\rho^{-1}(g)$}.

\medskip

{\bf Proof.} By the Tarski-Seidenberg theorem \cite[Theorem~2.2.1]{BCR} the image $\rho(\mathbb R^{3|V|})$ is a semi-algebraic subset of $\mathbb{R}^{|E|}$. Any semi-algebraic set admits a stratification by smooth manifolds \cite[Chapter 9]{BCR}. In particular, there exists an open set $R'$ of maximal dimension (and of full measure in $\rho(\mathbb R^{3|V|}$) which is a smooth submanifold of $\mathbb R^{|E|}$. The preimage $U:=\rho^{-1}(R')$ is an open subset of $R^{3|V|}$. Let $\Delta':=\Delta\cap U$. By stratification for algebraic varieties $\Delta'$ is a finite union $\cup_i \Delta_i$ of submanifolds  of positive codimension in $U$.

Restricting $\rho$ to $U$ one obtains a smooth map between two smooth manifolds $\rho\big|_U:U\to R'$. We can apply Sard's theorem to obtain a subset $R''\subset R'$ of full measure, such that any fiber $\rho^{-1}(g)$ for $g\in R''$ is smooth. 

Then for each $\Delta_i$ there is a subset $C_i\subset R''$ of full measure such that $\rho^{-1}(g)$ intersects $\Delta_i$ transversally for any $g\in C_i$. Indeed, the set of $g$, such that $\rho^{-1}(g)$ is not transversal to $\Delta_i$, has measure zero in $R''$ by Sard's theorem applied to $\rho\big |_{\Delta_i}$.

Transversality of intersection of $\rho^{-1}(g)$ with $\Delta_i$ implies that $\rho^{-1}(g)\cap \Delta_i$ is nowhere dense in $\rho^{-1}(g)$. In particular, it is nowhere dense in any connected component of $\rho^{-1}(g)$. Now $R:=\cap_i  C_i$ is a set which existence if claimed by the lemma.
$_\blacksquare$

\medskip

{\bf5.10.~Definition.} Each $r\in R$ gives rise to a polyhedral metric on $S$ (a sample polyhedron in the sense of Definition 2.12). We will call a sample polyhedron equipped with such a metric {\it generic}.

\medskip

{\bf5.11.~Corollary.} For a generic polyhedral metric on a surface $S$, the subset of its polyhedral realizations $q\in\mathbb E^S$, such that $q$ is boundary rigid, is open and dense in $\mathbb E^S$. We will denote it by $B_S\subset \mathbb E^S$.

\medskip

{\bf Proof.} Immediate from Proposition 5.7 and Lemma 5.9.$_\blacksquare$

\medskip

{\bf5.12.~Lemma.} {\sl Let $S$ be a generic sample polyhedron. Then the dimension of the image of $d\delta: \T_q\mathbb E^S\to \T_{q\circ \delta}\mathbb E^P$ at any $q\in B_S$ is  equal to $|F|+3$, where $|F|$ is the number of boundary edges. }

\medskip

{\bf Proof.} From exact sequence $\mathrm{ker}\,d\delta \to \T_q\mathbb E^S \to d\delta(\T_q\mathbb E^S)\to 0$ we obtain $\mathrm{dim}\,\delta(\T_q\mathbb E^S)=\mathrm{dim}\,\T_q\mathbb E^S-\mathrm{dim\, ker}\,\delta$. By Lemma 5.3 one has
$\mathrm{dim}\,\delta(\T_q\mathbb E^S)\ge |F|+3-\mathrm{dim\, ker}\,\delta$. On the other hand, the kernel of $d\delta$ at $q\in B_S$  is zero-dimensional by the definition of a boundary rigid polyhedron.
$_\blacksquare$

\medskip

Thus, taking the quotient by the full group of isometries ${\Isom}^+$ and combining with Theorem 4.3, we obtain the following

\medskip

{\bf5.13.~Theorem.}  {\sl For a generic polyhedral disk $S$, the image of the  map $d\delta: \T_q\mathbb E^S/{\Isom}^+\to \T_{q\circ\delta}\mathbb E^P/{\Isom}^+$ has dimension equal to $\frac{1}{2}\mathrm{dim}\,\T_{q\circ \delta}\mathbb E^P/{\Isom}^+=|F|-3$ on an open dense set. Thus $\delta(\mathbb E^S/{\Isom}^+)$ is Lagrangian semi-algebraic set (i.e. it contains an open dense Lagrangian submanifold).}
$_\blacksquare$ 

\bigskip

\noindent
\includegraphics[width=0.24\textwidth]{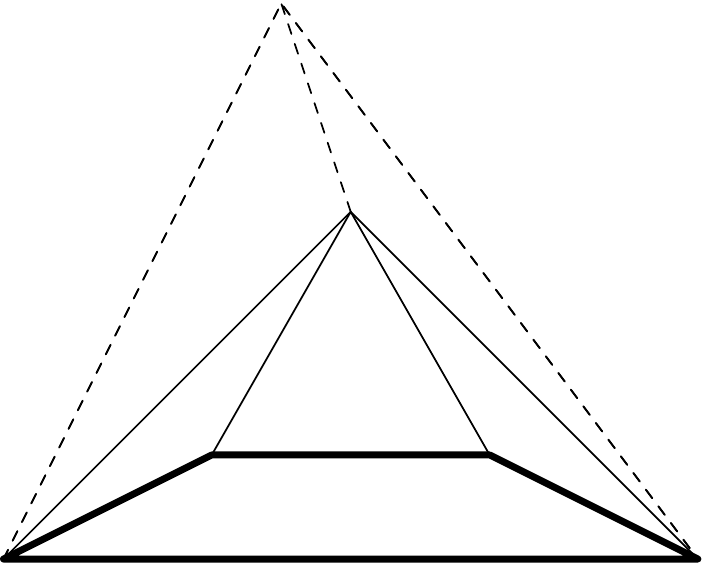}

\leftskip125pt

\vskip-100pt

\medskip
{\bf5.14.~Remark.} The genericity assumption is essential. As example on the left shows, a specific polyhedral surface with boundary can be rigid, that is, it may have a zero dimensional configuration space $\mathbb E^S/{\Isom}^+$ (see also \cite{AR}).  

The solid edges depict coplanar edges. The dashed edges intersect at a point not lying on the plane. The boundary is shown by bold lines.

\leftskip0pt

\vskip10pt

\noindent

\bigskip

\centerline{\bf5.~An application to a problem of R. Kenyon} %

\medskip

Richard Kenyon asked the following question on his web page \cite{Kenyon}:

\medskip

{\bf~Kenyon's problem.} {\sl Given a closed polygon $p$ in $\mathbb R^3$ composed of unit 
segments, is there an immersed polygonal surface
whose faces are equilateral triangles of edge length $1$, spanning p? }

\medskip

The conjecture was subsequentially resolved in the negative by Glazyrin and Pak \cite{GP}. Here we employ  the terminology of \cite{GP}, calling such polyhedral surfaces {\it domes} and polygons --- {\it integral curves}.
Ian Agol in a comment to a question of Mohammad Ghomi on Mathoverflow \cite{Agol}  sugested an alternative strategy to show that the statement of the conjecture is false for almost all integral curves. Namely, if one is able to show that the set of polygons swept out by the boundary of a given dome is ``isotropic'', and hence has measure zero, then the solution immediately follows from the fact that the set of spannable integral curves is a countable union of sets of measure zero.

\medskip

{\bf5.1~Remark.} The space $\mathbb E^P/{\Isom}^+$ possesses a natural measure, namely the symplectic measure assosiated to the symplectic form $\omega$. Since the  singularities are isolated \cite{KaM}, we simply ignore them in what follows.

\medskip

{\bf5.2~Lemma.} {\sl Given an orientable polygonal surface $S\subset \hat S$ with boundary $\partial D$ (as in Lemma 2.2) and corresponding sample polygon $P$, the subset $\delta(\mathbb E^P/{\Isom}^+)$ has measure zero}.

\medskip

{\bf Proof.} We are interested only in the set theoretical image of the space of polyhedra up to sets of measure zero. Thus we can take reduction of the scheme $\mathbb E^S/{\Isom}^+$ and consider its smooth locus which is open and of full measure. Now we have a smooth map of smooth manifolds, and by Theorem 4.3 the rank of its differential is at most half the dimension of the target manifold. Now it follows from the regularity (e.g. by the constant rank theorem) that the image of this map has measure zero in $\mathbb E^P/{\Isom}^+$.
$_\blacksquare$

\medskip

Now we are able to obtain a negative answer to the question of Kenyon in the case of (not necessarily immersed) orientable surfaces. This result extending to the non-orientable case was also proved in \cite{GP}.

\medskip

{\bf5.3~Theorem.} {\sl The set of integral curves in $\mathbb E^P/{\Isom}^+$ that are boundaries of orientable domes has measure zero.}

\medskip

{\bf Proof.}  The sample polyhedron of a dome is uniquely defined by its combinatorial structure, since all edges are assigned length $1$. There are countably many finite simplicial complexes. Thus the set of integral curves that can be spanned by domes is a countable union of sets of measure zero, hence has measure zero.
$_\blacksquare$

\medskip

Note that Glazyrin and Pak proved that this set is also dense.

\bigskip

\bibliographystyle{plain}
\bibliography{poly.bib}

\end{document}